\documentclass[10,twocolumn]{article}

\usepackage{fullpage}
\usepackage{amsmath}
\usepackage{cases}
\usepackage{amssymb}
\usepackage{amsthm}
\usepackage[all,cmtip,knot]{xy}
\usepackage{graphicx}
\usepackage{authblk}
\usepackage{epsfig}
\usepackage{epstopdf}
\usepackage{subfigure}
\usepackage{booktabs}
\usepackage{mathrsfs}
\usepackage{algorithmic}
\usepackage{algorithm}
\usepackage{setspace}
\usepackage{cite}

\usepackage{graphicx}
\usepackage{pnastwof}
\usepackage{epstopdf}
\usepackage{amssymb,amsfonts,amsmath}
\usepackage{subfigure}

\newcommand{\argmin}[1]{\underset{#1}{\mathrm{argmin}}}

\title{Sparse Dynamics for Partial Differential Equations}

\author{Hayden Schaeffer}
\author{Stanley Osher}
\author{Russel Caflisch}
\affil{\vspace{-0.05cm}University of California, Los Angeles, Department of Mathematics}

\author{\vspace{-0.25cm}Cory Hauck}
\affil{\vspace{-0.25cm} Oak Ridge National Laboratory}

\date{ \ }
\begin{document} \maketitle
\begin{abstract}
  We investigate the approximate dynamics of several differential equations when the solutions are restricted to a sparse subset of a given basis.  The restriction is enforced at every time step by simply applying soft thresholding to the coefficients of the basis approximation. By reducing or compressing the information needed to represent the solution at every step, only the essential dynamics are represented. In many cases, there are natural bases derived from the differential equations which promote sparsity. We find that our method successfully reduces the dynamics of convection equations, diffusion equations, weak shocks, and vorticity equations with high frequency source terms.
\end{abstract}

\maketitle

\section{Introduction}
In this work, we investigate the approximate dynamics of various PDEs whose solutions exhibit behaviors on multiple spatial scales.  These scales may interact with one another in a non-linear manner as they evolve. Many physical equations contain multiscale (as well as multiphysics) phenomena,  such as the homogenization problems from material science and chemistry and multiscale systems in biology, computational electrodynamics, fluid dynamics, and atmospheric and oceanic sciences. In some cases, the physical laws used in the model can range from molecular dynamics on the fine scale to classical mechanics on the large scale. In other cases, the equations themselves contain high wave number oscillations that separate into discrete scales, on top of the smooth underlying behavior of the system.

The main source of difficulty in multiscale computation  is that accurate simulation of the  system requires   all phenomena to be fully resolved. The smaller spatial scales influence the global solutions, thus they cannot be ignored in the numerical computation. In some cases it is possible to derive an analytic equation for the  effect of small scales on the solution (see \cite{arnold_mathematical_1989,papanicolau_asymptotic_1978}). In practice, however, it may not be possible to derive a simple expression that represents the fine scale behavior.  Many problem dependent methods have been proposed in the literature, while a few provide a general methodology for modeling the macroscopic and microscopic processes that yield multiscale models. For example, some general methods include the heterogenous multiscale method  \cite{weinan2003heterogeneous}, the equation-free method \cite{kevrekidis2003equation}, multiscale methods for elliptic problems \cite{nolen2008framework}, and the sparse transform method \cite{daubechies_sparse_2007}. For an overview of general multiscale approaches see \cite{e_multiscale_2003}. A key difference between our method and the methods in \cite{kevrekidis2003equation, nolen2008framework, weinan2003heterogeneous} is that we are directly resolving all of the significant scales in the solution. By contrast,  the methods of \cite{kevrekidis2003equation, nolen2008framework,weinan2003heterogeneous} directly resolve only the coarse scales of the solution, and they separately ``reconstruct" the fine scale solution
(as well as its effect on the coarse scales).

From the perspective of mathematics, multiscale methods began with representation of a function using a global basis, such as  Taylor series or Fourier series. More sophisticated bases have appeared; for example, any one of the many wavelet bases used in imaging and computational physics.  The key to the basis approximation is that each basis element represents behavior on a specific  scale, therefore the coefficients of the basis provide complete information about the underlying function. This is also the principle behind multiresolution and decomposition methods. 

As the methods of multiscale and multiphysics modeling developed over the past few decades, so did corresponding methods in imaging and information science. One of the fundamental ideas in imaging is that of sparsity. Sparse data representation is used throughout imaging from compression to reconstruction.
Early advances in sparse techniques were, e.g., in  \cite{candes2006, donoho2006compressed},   which presented a convex minimization approach to the computationally challenging sparse basis pursuit problem. Many models have been proposed which use sparsity to produce both more efficient numerical methods and better quality solutions. Some applications of sparsity to imaging include: compressive sensing, reconstruction of images from sparse data \cite{zhou2009non,schaeffer2011low} ,  and recovery of images using sparse regularization \cite{goldstein2009split,rudin1992nonlinear}. The underlying principle of sparsity is that images can be approximately represented by a small number of terms with respect to some basis. Inducing sparsity, creating effective bases, and developing efficient computational algorithms have been intensely active fields in information science. 

For imaging and information science, one of the reasons for the success of sparse methods  is their ability to resolve drastically different phenomena with a small amount of information. This is also a principal goal of multiscale modeling. In this work, we transfer  sparsity methodology, which was developed for  information science, to multiscale nonlinear differential equations and show that it can be an effective tool for accurately computing solutions using less information.

In particular, we propose solving PDEs with the constraint that the approximate solution resides on a sparse subspace of a basis. In this way, the complexity of the method will depend on the number of basis terms retained and (\textit{nearly}) independent of the grid size. In the following sections we will discuss the general problem and the optimization method used to induce sparsity in the solution. Also, the general numerical method will be explained as well as results of numerical experiments. The method is tested on an advection equation with oscillatory velocity, a parabolic equation with oscillatory coefficients, a conservation law with oscillatory diffusion, and the vorticity equations with high frequency source terms. We conclude with a discussion on the proposed work and implications.

\section{Sparsity}

\subsection{Problem Statement}
In general, the problem can be stated as follows. Assuming $x \in \mathbb{R}^n$ and $t>0$, let $u(x,t): \Omega \rightarrow \mathbb{R}$  be the approximate solution of 
\begin{equation}
\frac{\partial u}{\partial t}= F(u) 
\label{eq:Sparse}
\end{equation}

\noindent subject to the constraint: $u(x,t) =\sum_j \widehat{u}_j(t) \phi_j(x)$, where the number of non-zero $c_j$ at a given time step is sparse. The operator $F(\cdot)$ can be non-linear, non-local, and dependent on the derivative of $u$. The basis terms $\phi_j$ are assumed to exist on separate scales, which is true of most bases (for example, Legendre polynomials, Fourier, Wavelets, \textit{etc.}). In this way, the basis terms represent different global behaviors.  The method involves two steps: evolve the PDE forward in time and project the updated solution onto a sparse subset. We first address sparsity induction through soft thresholding.

\subsection{Sparsity via Optimization: Soft Thresholding}
At a given time step, the problem of projecting the updated solution onto a sparse subset is equivalent to fitting a solution $u^{n}$ with corresponding coefficients $\{\widehat{u}^n_j\}$ at the $n$th time step to a solution $u$ whose corresponding coefficients $\{\widehat{u}_j\}$ are sparse. This can be written as a constrained least squares fit as follows.
\begin{eqnarray}
 & & \hspace{1cm} \min_u \ ||u-u^n||^2_{L^2} \ \ \ \ \ \text{s.t.}\label{eq:LeastSquares} \\ 
& & u =\sum_j \widehat{u}_j \phi_j \ \ \&\ \  \{ \widehat{u}_j \} \ \ \text{is sparse} \nonumber
\end{eqnarray}

\noindent Expanding with respect to the basis and assuming that the basis is orthonormal, this constrained optimization problem is related to the following unconstrained problem:

\begin{eqnarray}
 \text{(L0)} \hspace{0.4cm}\widehat{u}=\argmin{\widehat{u}} \ \ \lambda ||\widehat{u}||_{0}+\frac{1}{2} ||\widehat{u}-\widehat{u}^n||^2_{L^2} 
\label{eq:L0}
\end{eqnarray}

\noindent where $\widehat{u}$ is the vector of coefficients. The ``norm'' $||\cdot||_0$ is the number of non-zero coefficients in equation \ref{eq:LeastSquares}. This makes equation \ref{eq:L0} both non-convex and difficult to solve. By replacing the $L^0$ norm by the $L^1$ norm, we get the following convex relaxation of equation \ref{eq:L0}:

\begin{eqnarray}
\text{(L1)} \hspace{0.4cm} \widehat{u}=\argmin{\widehat{u}} \ \ \lambda ||\widehat{u}||_{1}+\frac{1}{2}  ||\widehat{u}-\widehat{u}^n||^2_{L^2} 
\label{eq:L1}
\end{eqnarray}

\noindent Note that since $\widehat{u} \in \mathbb{C}$, the $L^1$ norm is 
$||\widehat{u}||_{L^1} := \sum_j |\widehat{u}_j|$, where $|\widehat{u}_j|=\sqrt{ \widehat{u}_j \ \overline{\widehat{u}}_j}$. The solution of equation \ref{eq:L1} is given by the following equation: 
\begin{eqnarray}
\widehat{u}_j &=& \mathcal{S}_{\lambda} \left( \widehat{u}_j^n\right) \nonumber \\
&=& \mathrm{max} \left( |\widehat{u}_j^n|- \lambda, 0 \right) \ \frac{\widehat{u}_j^n}{|\widehat{u}_j^n|}
\label{eq:Shrink1}
\end{eqnarray}

In general, this can be computed for a non-orthonormal basis, which is equivalent to a basis pursuit problem with the $L^1$ norm as a sparse regularizer. In that case, the solution must be found by an iterative method rather than the simple shrinkage provided here as an example. The resulting minimizer $\widehat{u}$ is a proximal solution which lies on a sparse subset of the original coefficient domain \cite{cai2009split}. This can be used to show that the solutions form a contraction map in the $L^2$ norm. Alternatively, we can simply apply the soft thresholding on the coefficients directly in order to induce sparsity in this way. 

\section{Numerical Method}
Assuming $u(x,t)$ is periodic in the domain $\Omega \subset \mathbb{R}^n$, one natural basis is the Fourier basis, whose coefficients are the Fourier transform of $u(x,t)$. This is appropriate for the examples shown here. For the rest of this work, we will use the Fourier basis; however, the overall methodology presented here is independent of the corresponding basis.

Taking the Fourier transform of the PDE from equation \ref{eq:Sparse} and discretizing the resulting differential equation in time yields a multistep scheme. Since our method does not depend on the choice of numerical updating, we can assume that the scheme takes the following form:
\begin{equation}
\widehat{v}= Q \left(\widehat{u}^{n-q}, ..., \widehat{u}^{n} \right)
\end{equation}
The updated solution $\widehat{v}$ may be sparse depending on both the PDE and the update operator $Q$, but in general will have non-trivial values everywhere depending on the approximation and implementation. The auxiliary variable $v$ is projected onto a sparse subspace by the shrinkage operator:
\begin{equation}
\widehat{u}^{n+1}= \mathcal{S}_{\lambda} \left(\widehat{v} \right)
\end{equation}

\noindent Altogether, the update in the spatial domain is simply:

\begin{equation}
u^{n+1}=  \sum_j \mathcal{S}_{\lambda} \left(Q \left(\widehat{u}^{n-q}, ..., \widehat{u}^{n} \right) \right) \phi_j
\end{equation}
Unlike traditional projections, this is non-linear and adaptive. Rather than sorting the coefficients and retaining a fixed number of large amplitude terms or keeping terms whose wavenumbers are below some cutoff, the shrinkage allows the number and choice of non-zero coefficients to evolve over time. Also, this is not the same as hard thresholding the solution at every step (\textit{i.e.} keeping only the terms larger than a fixed value) since the coefficients that remain have decreased their magnitude by $\lambda$. Most importantly, the projection does not favor any particular part of the spectrum; instead the amplitude of the coefficient determines if it remains. In terms of the Fourier basis, the importance is placed on the amplitude not the wavenumber.

For general convergence, as long as $\lambda=C dt^{p}$ for $p$ larger than the accuracy of the scheme used to update the variable in time, then the shrinkage operation does not change the spatial accuracy of the original method and the method will still converge as $dt \rightarrow 0$.  For example, discretize using the forward Euler method, and then expand the shrinkage operator to get
\begin{equation}
\widehat{u}^{n+1}=\widehat{u}^{n}+ dt \ \widehat{F(u^n)}+\mathcal{O} \left( \lambda \right).
\end{equation}
Therefore, to have convergence as $dt \rightarrow 0$, the shrinkage parameter must be $\lambda=C dt^{1+\alpha}$. In general, the shrinkage operator is non-expansive in each coefficient, hence non-expansive in coefficient norms. This may help with obtaining a general convergence result.

\section{Numerical Results}

In this section, we discuss the application of the proposed sparse method to several equations with different numerical schemes.
\subsection{Convection}
The convection equation we consider is the following:
\begin{equation}
\partial_t u =a(x) \partial_x u
\label{eq:Hyper}
\end{equation}

\noindent where the coefficient $a(x)$ is highly oscillatory.
 
\begin{figure*}[t]
\begin{center}
\subfigure[ True (black) and Sparse (blue) Solution in $x$ space]{
\includegraphics[width=2.5in]{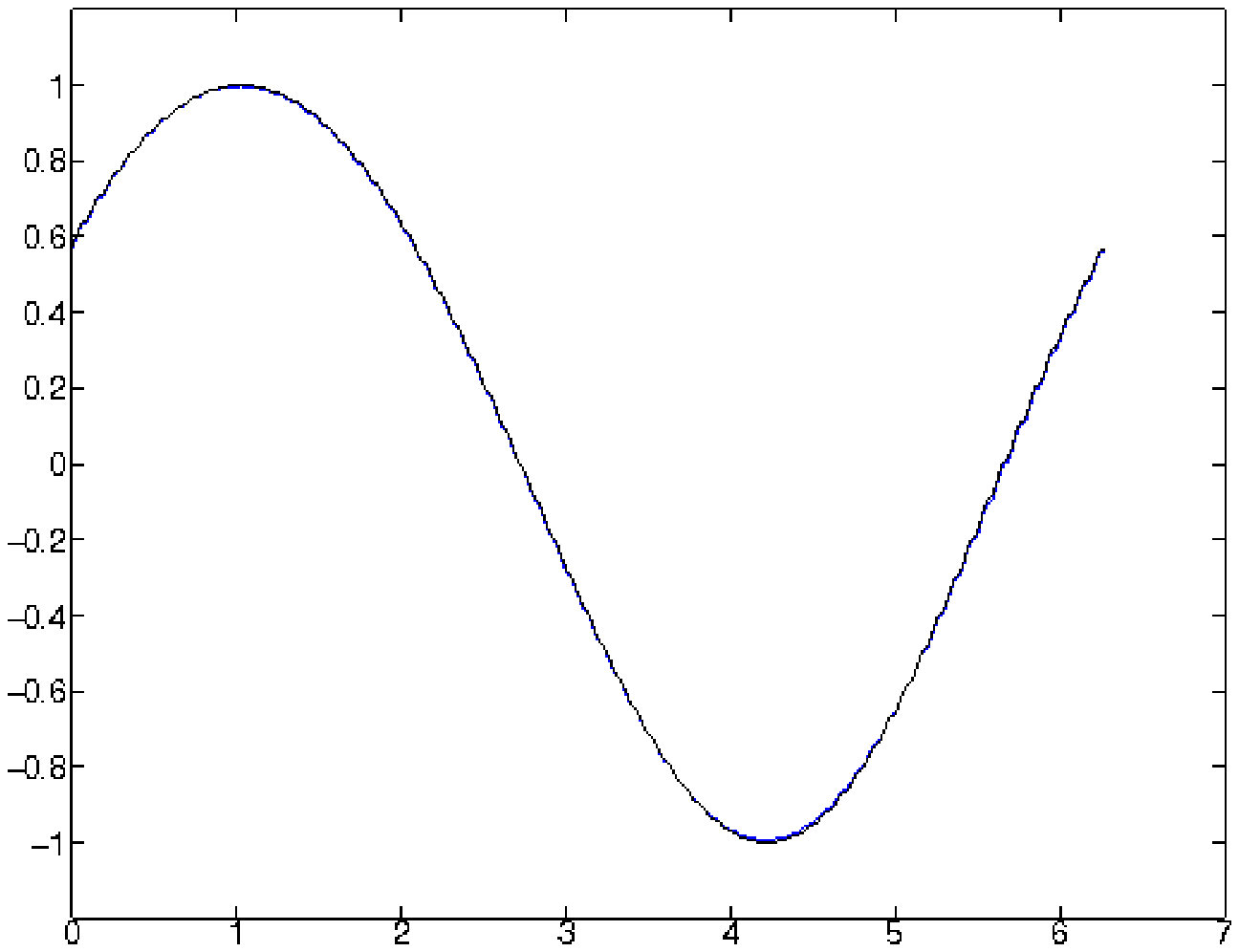}
} \hspace{1cm}
\subfigure[True (black) and Sparse (blue 'x') Solution in $x$ space, zoomed in]{
\includegraphics[width=2.5in]{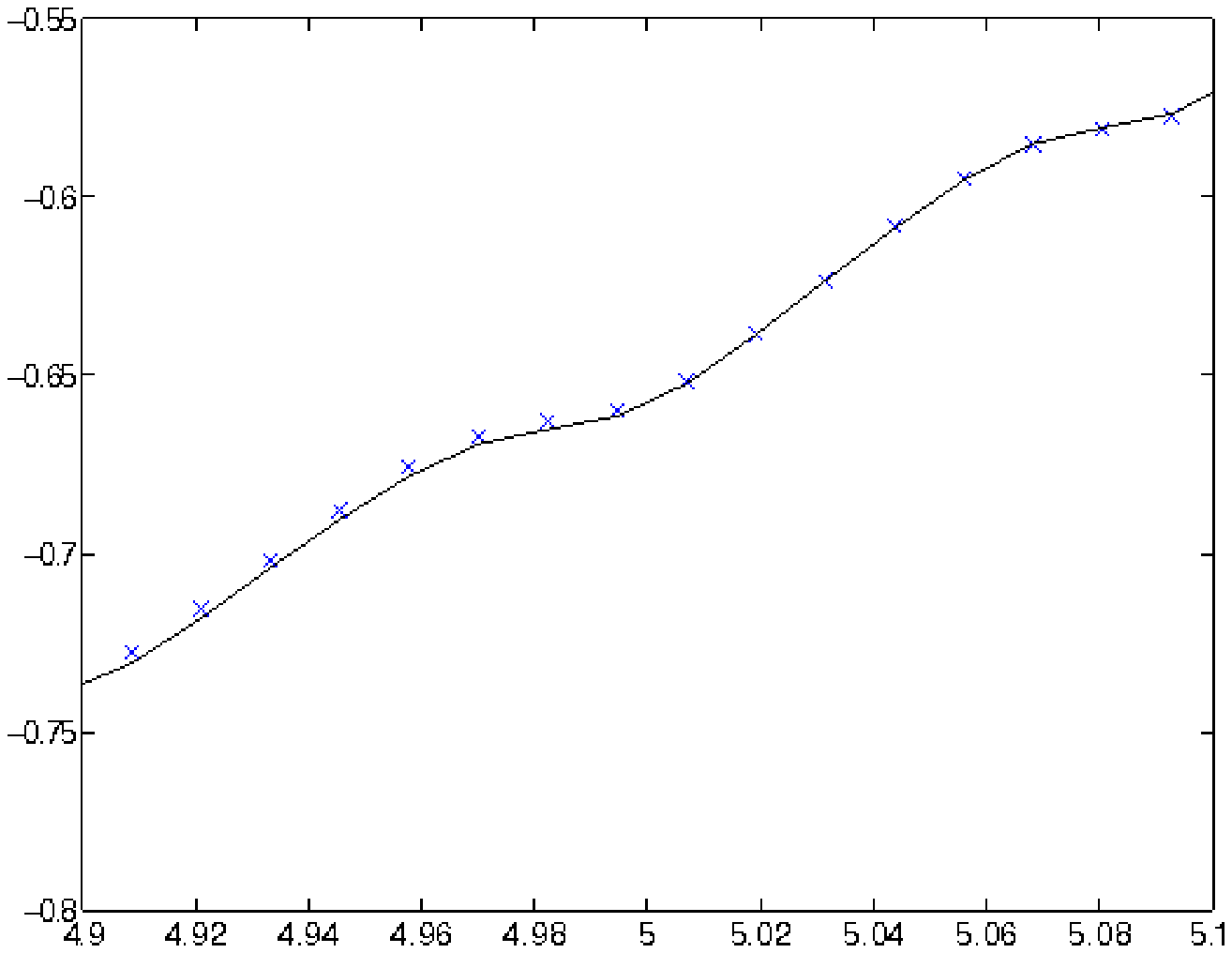}
}\hspace{-0cm}
\subfigure[Sparse (blue) Solution and low frequency (red) Solution in $x$ space]{
\includegraphics[width=2.5in]{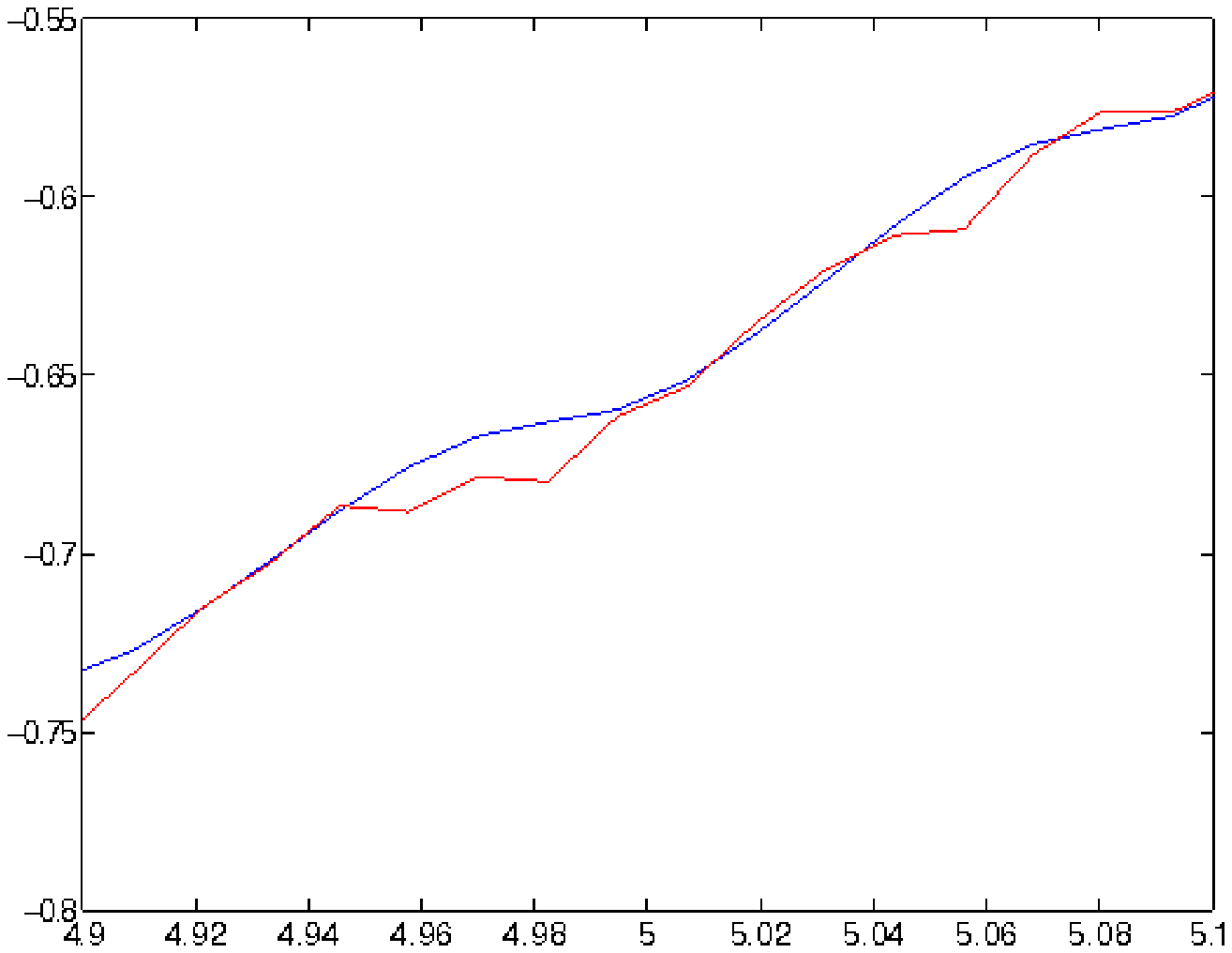}
}\hspace{1cm}
\subfigure[True (black) Solution and Sparse (blue) Solution in $\widehat{u}$ domain]{
\includegraphics[width=2.5in]{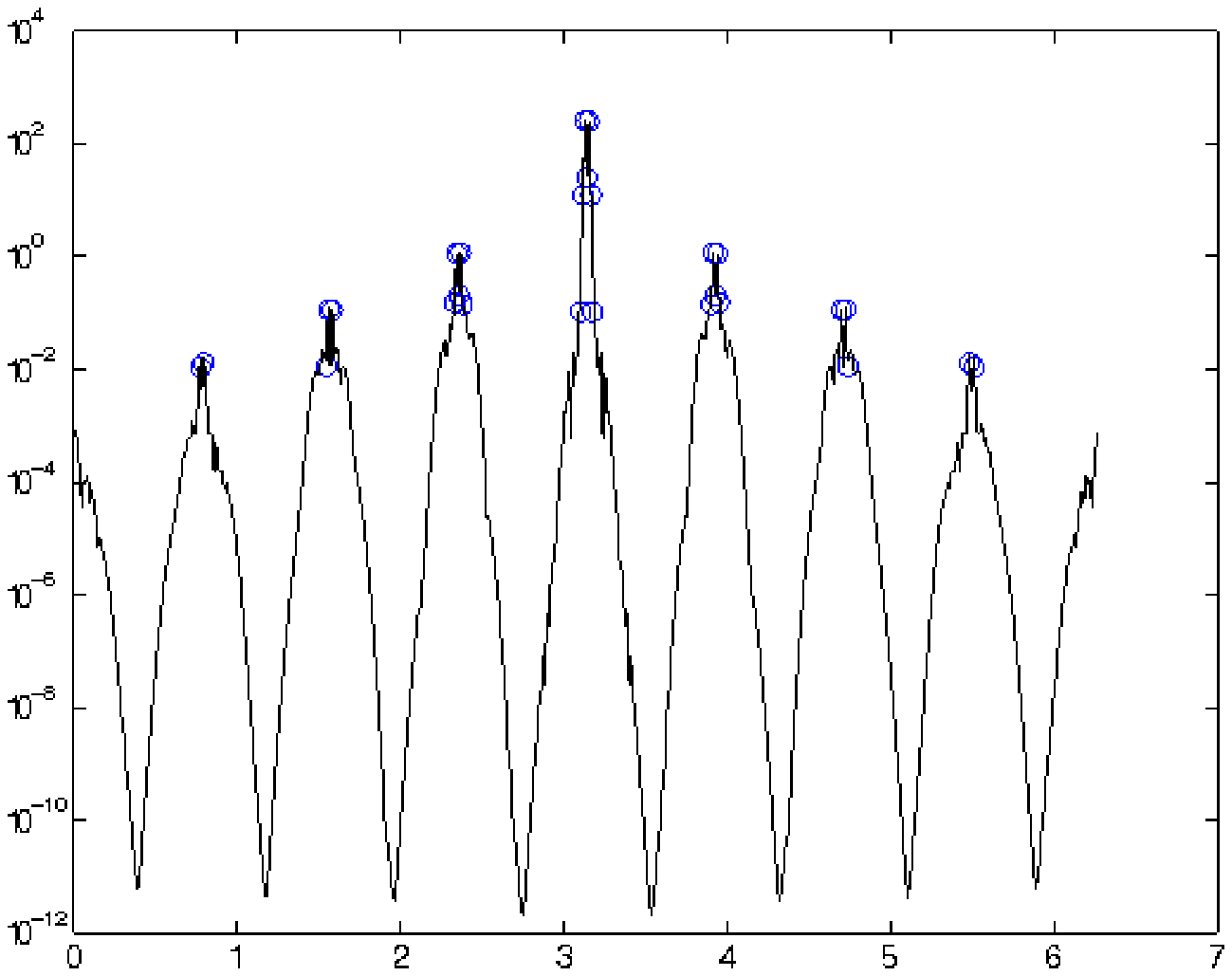}
}
\caption{Convection with highly oscillatory coefficients. The solution is shown in $x$ space, zoomed in version in $x$ space, and in the $\widehat{u}$ space. The solutions are shown on a 512 grid, with $dt=2e-3$, $dx=1.23e-2$, and $\lambda=5e-03$.}
\label{fig:Hyper}
  \end{center}
\end{figure*}

Let $k$ be the wavenumber and use spectral Leap Frog as the updating for equation \ref{eq:Hyper} to obtain

\begin{equation}
 \widehat{v} ^{n+1}=\widehat{u}^{n-1}+ \ 2 \ dt \ \ \widehat{a} \ast \left( \ i \ k  \ \widehat{u}^{n} \right)
 \label{eq:HyperUpdate}
\end{equation}
in which $\ast$ is the convolution operator over frequency.

The time step is $\mathcal{O}(dx)$ to preserve the stability condition in equation \ref{eq:Hyper}. In Figure \ref{fig:Hyper}, the coefficient is chosen as follows.
\begin{equation*}
a(x)={\frac{1}{4}} \exp\left({\frac{0.6+0.2 \cos(x)}{1+0.7 \sin(64x)}}\right)
\end{equation*}

\noindent This choice of $a(x)$ exhibits both fast and slow modes, but the particular structure is not directly needed.

Figure  \ref{fig:Hyper} illustrates the performance of the sparse solution method on this example by comparison of the sparse solution, the true solution produced using a standard fully resolved method, and a ``low frequency solution" produced by solving equation \ref{eq:HyperUpdate} for wavenumbers $k$ in the interval $|k|\leq K$ in which $K$ is the number of modes in sparse solution.
 In Figure \ref{fig:Hyper} (a-b) the sparse solution produced by our method and the true solution at $t=1$ are plotted in the spatial domain at a given time. In Figure \ref{fig:Hyper} (c) the sparse and low frequency solutions are displayed. The low frequency solution contains parasitic modes common to the Leap frog scheme. In Figure \ref{fig:Hyper} (d) the sparse and true spectra are plotted. The sparse spectrum captures the largest amplitude coefficients throughout the domain. In fact, out of the 512 coefficients used in the true solution, only 27 are retained in the sparse one (about $5.3\%$). 

\subsection{Parabolic}
The parabolic equation we consider is the following:
\begin{equation}
\partial_t u = \partial_x \left( a(x) \partial_x u \right)
\label{eq:Parabolic}
\end{equation}

\noindent where the diffusion coefficient $a(x)$ is highly oscillatory. The coefficient is assumed to be bounded, \textit{i.e.} $A_M \geq a(x) \geq A_m > 0$. This is also related to the elliptic case $\partial_x \left( a(x) \partial_x u \right)=f$, since an elliptic equation can be solved by taking a parabolic scheme to steady state. Alternatively, the corresponding parabolic scheme can be iterated forward for a small number of time steps in order to find a partial solution to the elliptic problem. Then  by using the partial solution, the locations of the non-zero coefficients can be extracted and the elliptic problem can be solved by a Galerkin method on these coefficients \cite{daubechies_sparse_2007}.
 
\begin{figure*}[t]
\begin{center}
\subfigure[ True (black) and Sparse (blue) Solution in $x$ space]{
\includegraphics[width=2.6in]{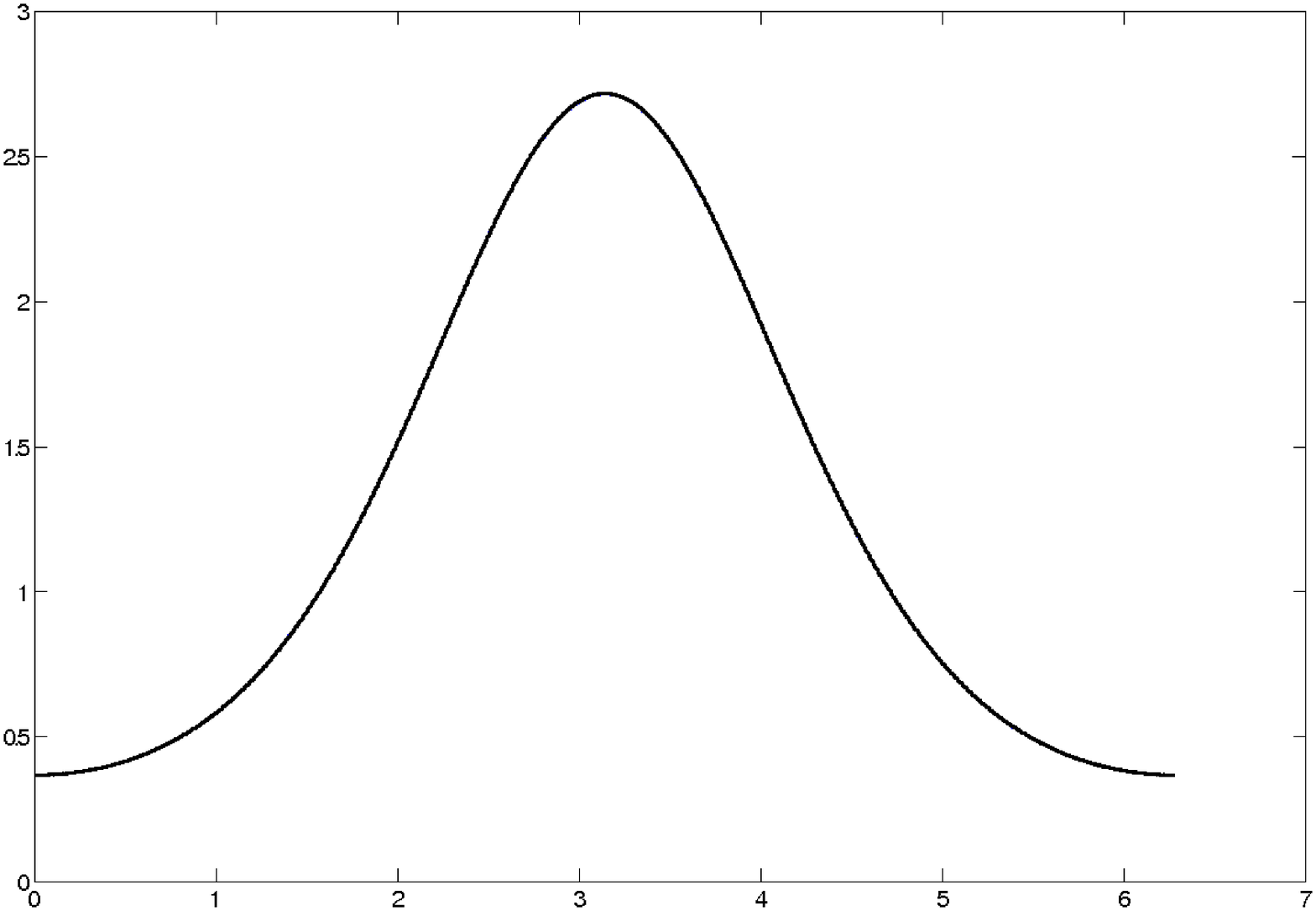}
} \hspace{1cm}
\subfigure[True (black) and Sparse (blue 'x') Solution in $x$ space, zoomed in]{
\includegraphics[width=2.6in]{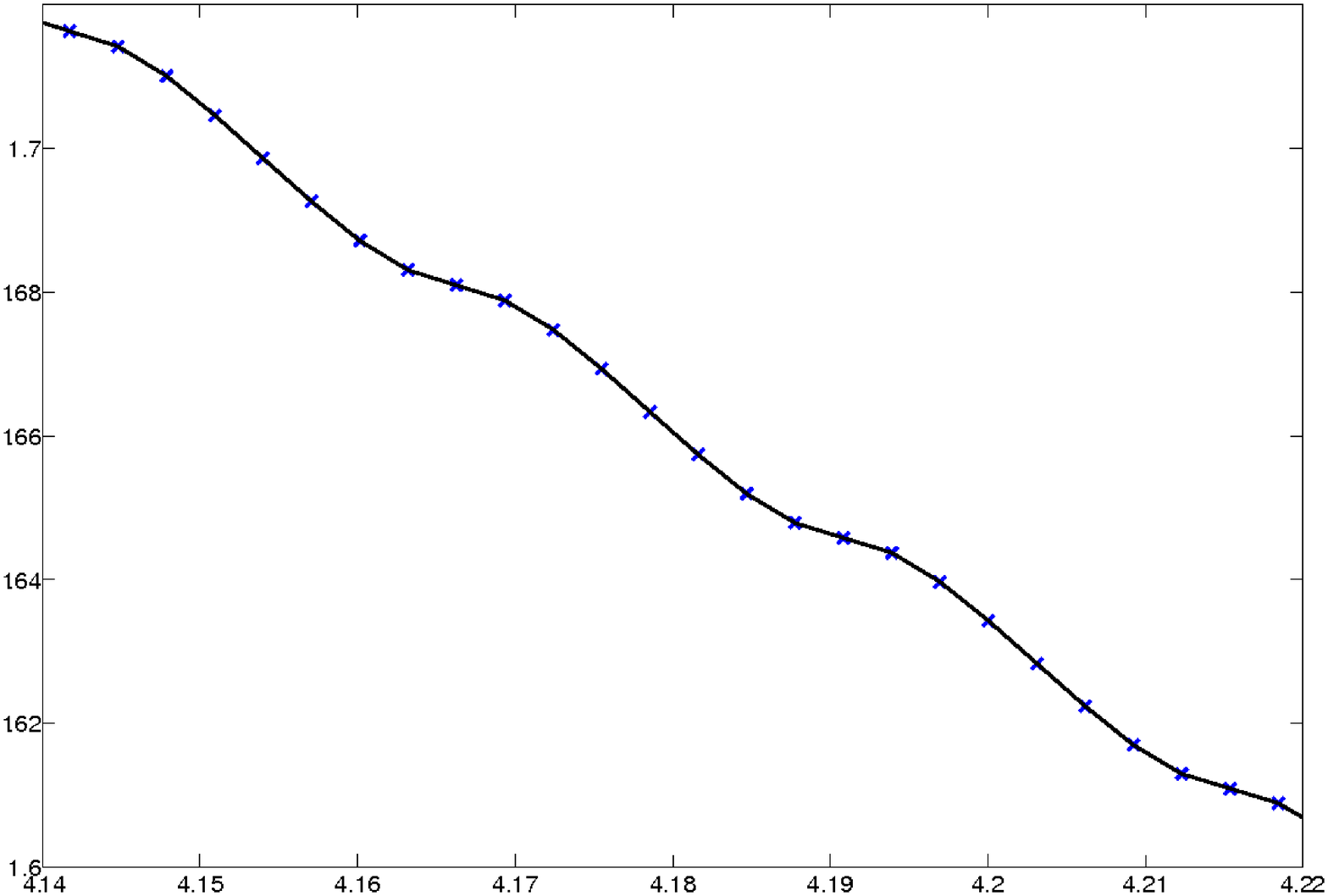}
}\hspace{-0cm}
\subfigure[True (black) Solution and Sparse (blue) Solution in $\widehat{u}$ domain]{
\includegraphics[width=2.6in]{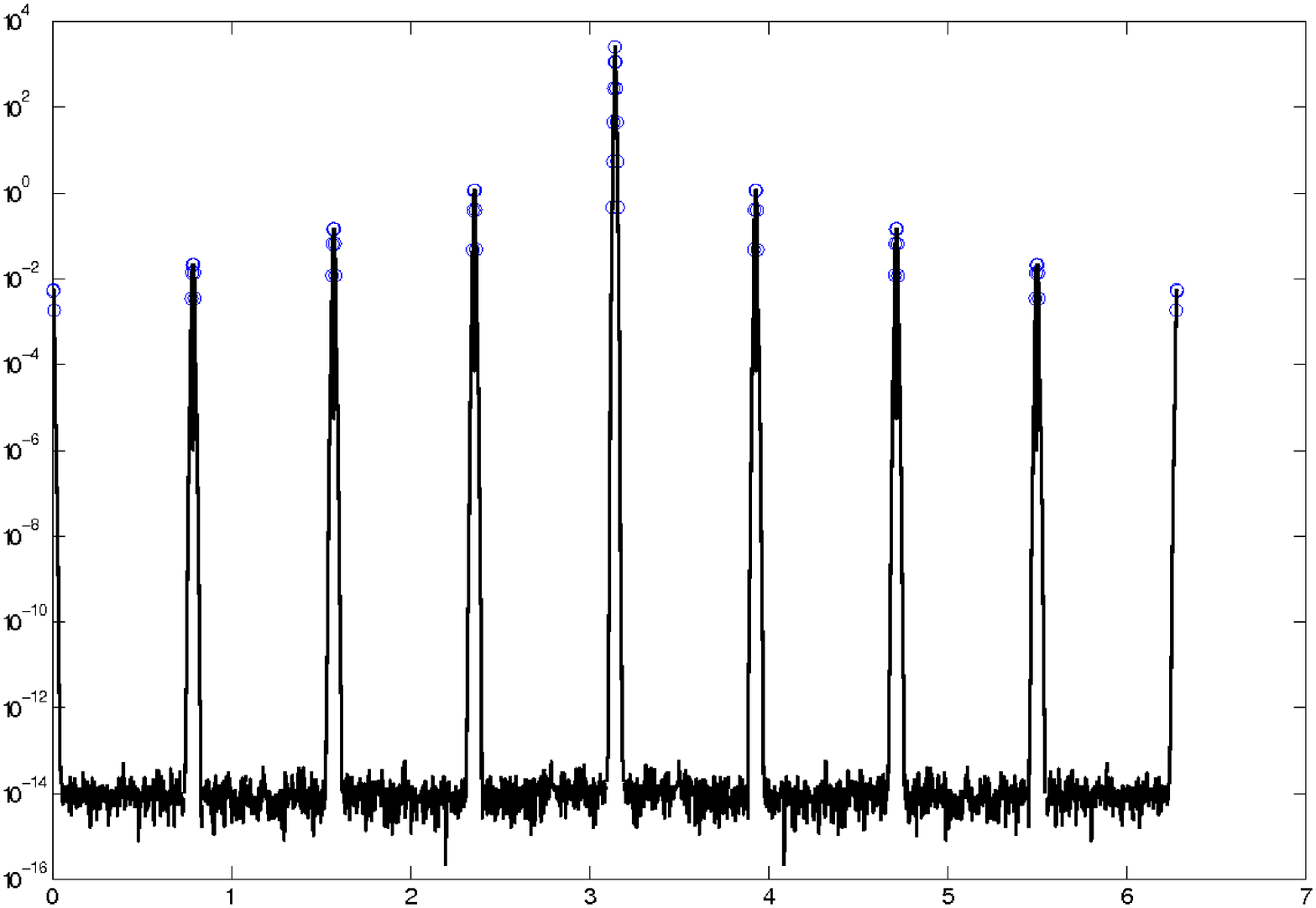}
}\hspace{1cm}
\subfigure[Diffusion coefficient in $x$ space]{
\includegraphics[width=2.6in]{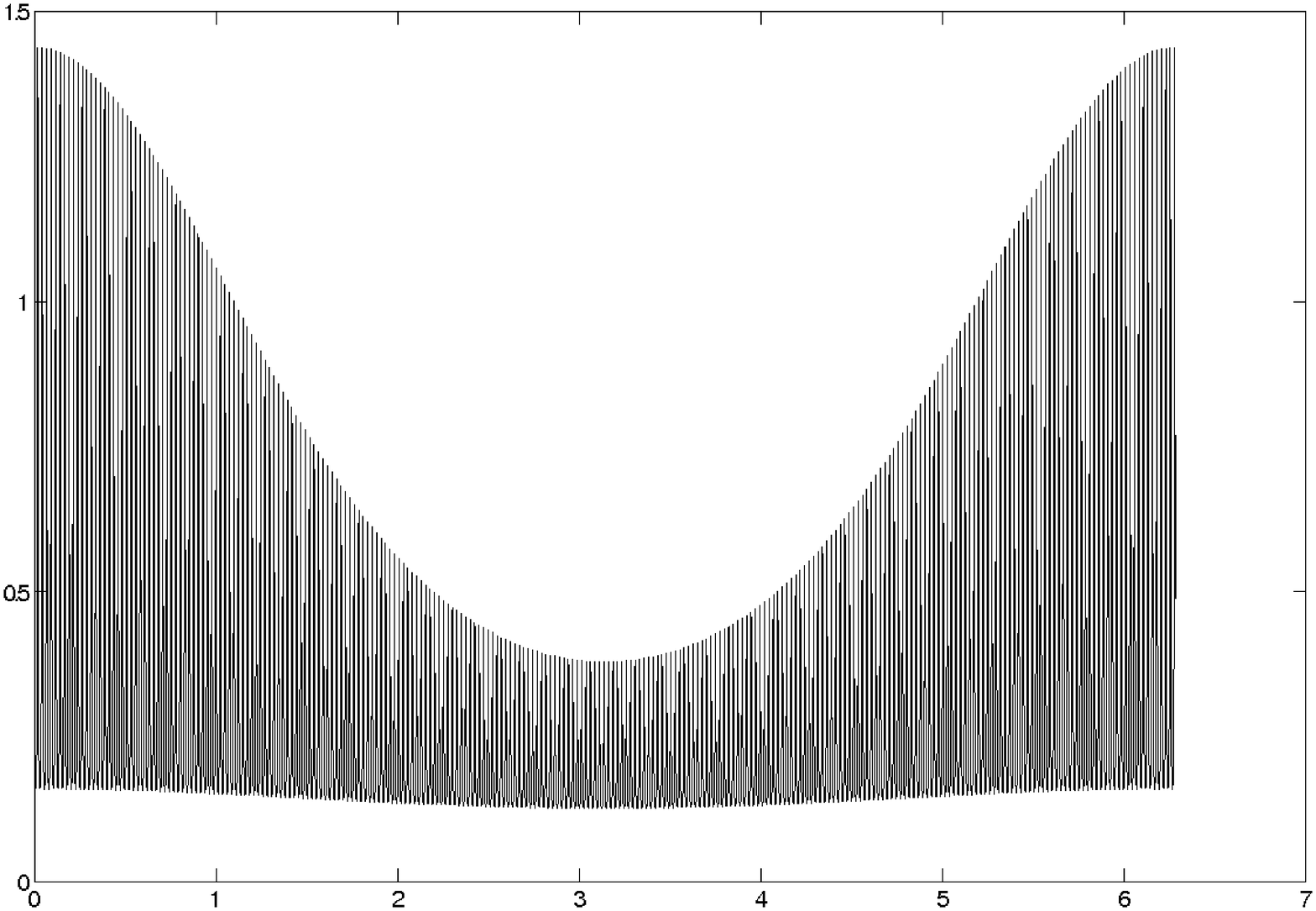}
}
\caption{Parabolic diffusion with highly oscillatory coefficients. The solution is shown in $x$ space, zoomed in version in $x$ space, and in the $\widehat{u}$ space. The solutions are shown on a 2048 grid, with $dt=1.5e-8$, $dx=3.1e-3$, and $\lambda=2.5e-06$.}
\label{fig:Parabolic}
  \end{center}
\end{figure*}

The updated scheme we use for equation \ref{eq:Parabolic} is forward Euler.

\begin{equation}
 \widehat{v} ^{n+1}=\widehat{u}^{n}+ \ dt \ i \ k  \ \ \widehat{a} \ast \left( \ i \ k  \ \widehat{u}^{n} \right)
 \label{eq:ParabolicUpdate}
\end{equation}

The time step is $\mathcal{O}(dx^2)$ to preserve the stability condition as well as the highly oscillatory nature of the coefficient $a(x)$ in equation \ref{eq:Parabolic}. In Figure \ref{fig:Parabolic}, the coefficient is chosen as follows.
\begin{equation*}
a(x)={\frac{1}{10}} \exp\left({\frac{0.6+0.2 \cos(x)}{1+0.7 \sin(256x)}}\right)
\end{equation*}

In Figure \ref{fig:Parabolic} (d) the highly oscillatory diffusion coefficient $a(x)$ is plotted in space. In Figure \ref{fig:Parabolic} (a-b) the sparse solution produced by our method and the true solution  at $t=1$ are plotted in the spatial domain at a given time and are nearly indistinguishable. The high frequency information is near the scale of the grid size, which can be seen in the zoomed in plot. In Figure \ref{fig:Parabolic} (c) the true and sparse spectra are displayed. The sparse spectrum captures the largest coefficients throughout the domain and not just the low wavenumbers. In fact, out of the 2048 coefficients used in the true solution, only 53 are retained in the sparse one (about $2.6\%$). In time, this number of non-zero coefficients as well as the identities of the non-zero coefficients will change in order to capture various behaviors.

\subsection{Viscous Burgers}

To investigate the sparse dynamics of conservative laws with diffusion, we use the viscous Burgers type equation.
\begin{equation}
\partial_t u + \frac{1}{2} \partial_x \left( u^2 \right) = \partial_x \left( a(x) \partial_x u \right)
\label{eq:Burgers}
\end{equation}

\noindent The LHS of equation \ref{eq:Burgers} is the standard Burgers advection term and the RHS is diffusion related to equation \ref{eq:Parabolic}. The equation exhibits three separate phenomenon: 1. Smooth large scale behavior from the diffusion, 2. Small scale oscillations induced from the high frequencies in the coefficient $a(x)$, and 3. Non-linear interactions between frequencies from the advection term. The update scheme in time is the standard total variational diminishing Runge-Kutta 2.
\begin{eqnarray*}
  \widehat{u}_{1}&=&\widehat{u}^{n}+ \ dt \ i \ k \left( \ \ \widehat{a} \ast \left( \ i \ k  \ \widehat{u}^{n} \right)- \widehat{F(u^{n})} \right)\\
  v ^{n+1}&=&\frac{1}{2} \left( \widehat{u}^{n}+ \widehat{u}_{1} \right)+ \ \frac{dt}{2} \ i \ k  \left( \ \widehat{a} \ast \left( \ i \ k  \ \widehat{u}_{1} \right)- \widehat{F(u_{1})} \right)
\end{eqnarray*}

\noindent where $F(u)=\frac{1}{2}u^2$. As before, we have the stability condition $dt$ is $\mathcal{O}(dx^2)$.

\begin{figure*}[tb]
\begin{center}
\subfigure[ True (black), Sparse (blue), and Low Frequency (red) Solution in $x$ space]{
\includegraphics[width=2.35in]{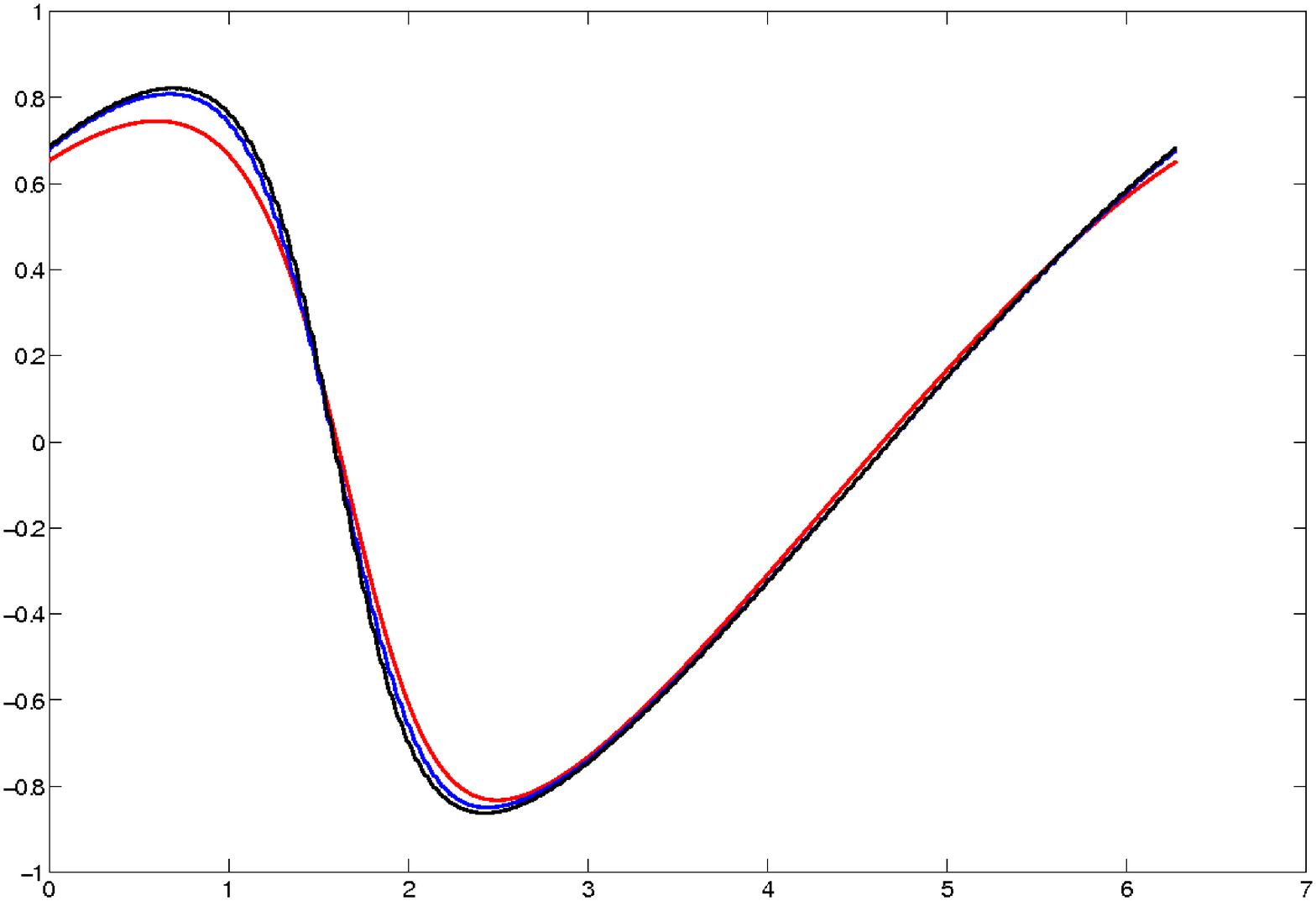}
} \hspace{1cm}
\subfigure[True (black), Sparse (blue), and Low Frequency (red) Solution in $x$ space, zoomed in]{
\includegraphics[width=2.35in]{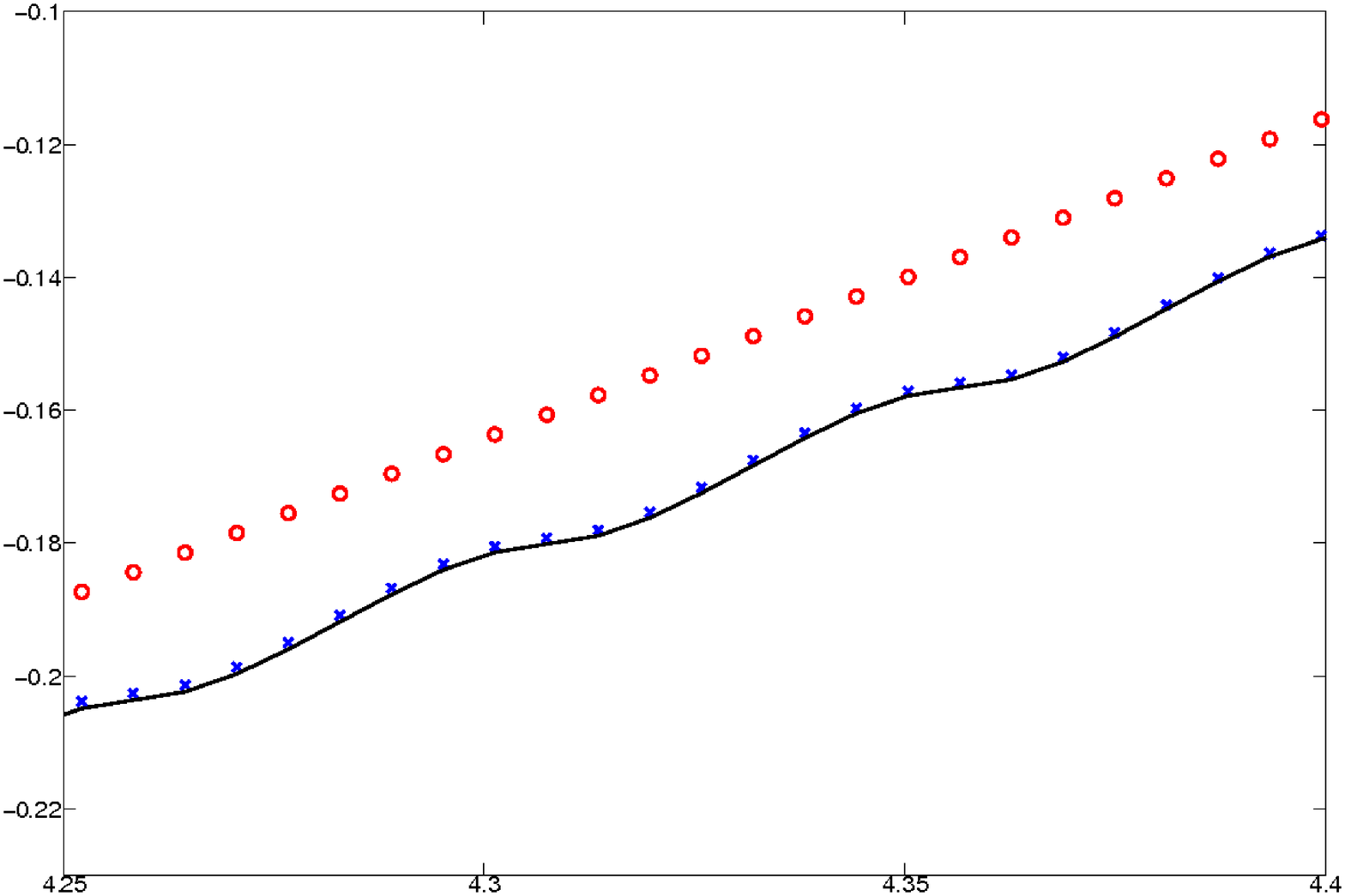}
}\hspace{-0cm}
\subfigure[True (black) and Sparse (blue) solutions in $\widehat{u}$ domain]{
\includegraphics[width=2.35in]{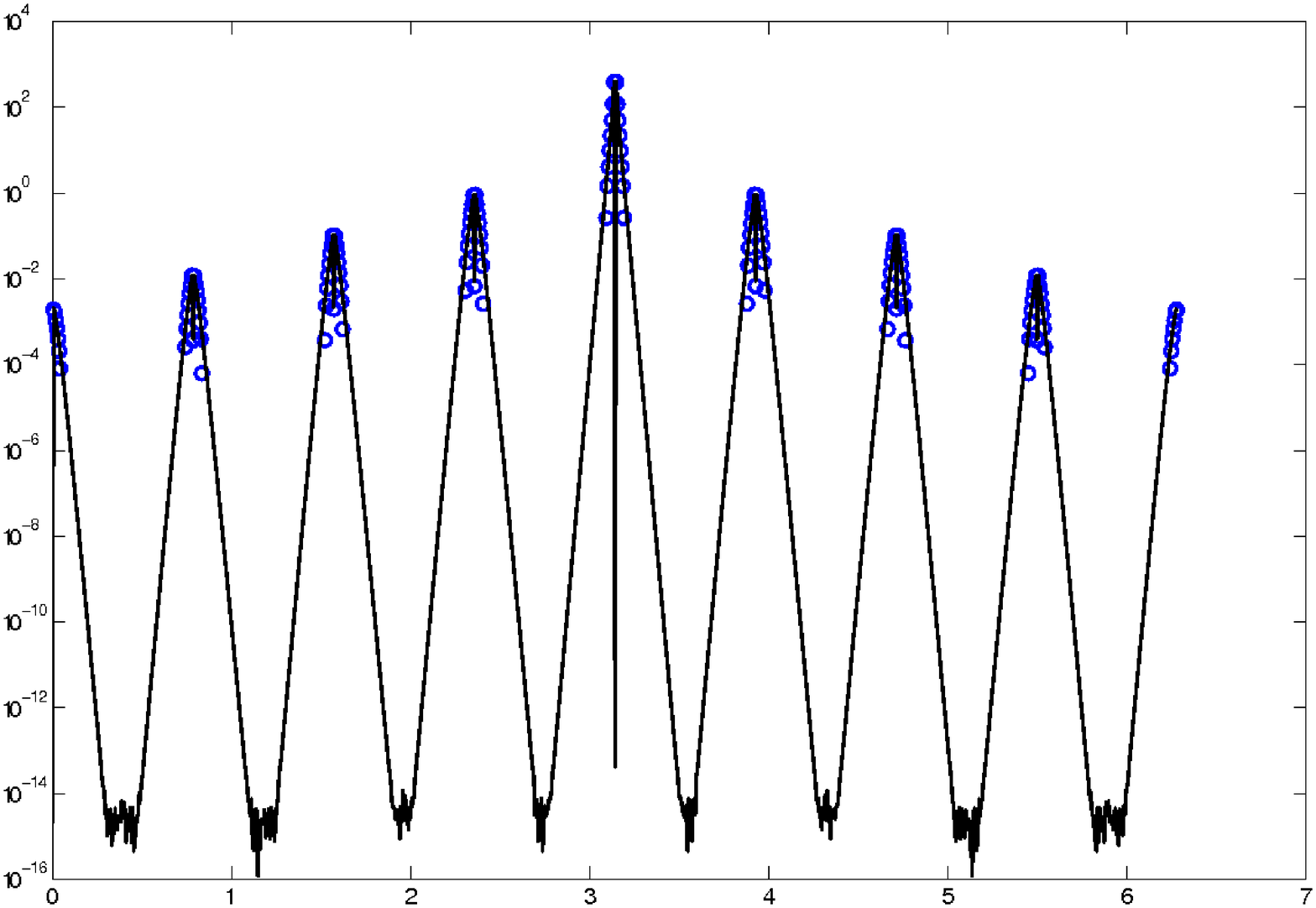}
}
\hspace{1cm}
\subfigure[True (black) and Low Frequency (red) solutions in $\widehat{u}$ domain]{
\includegraphics[width=2.35in]{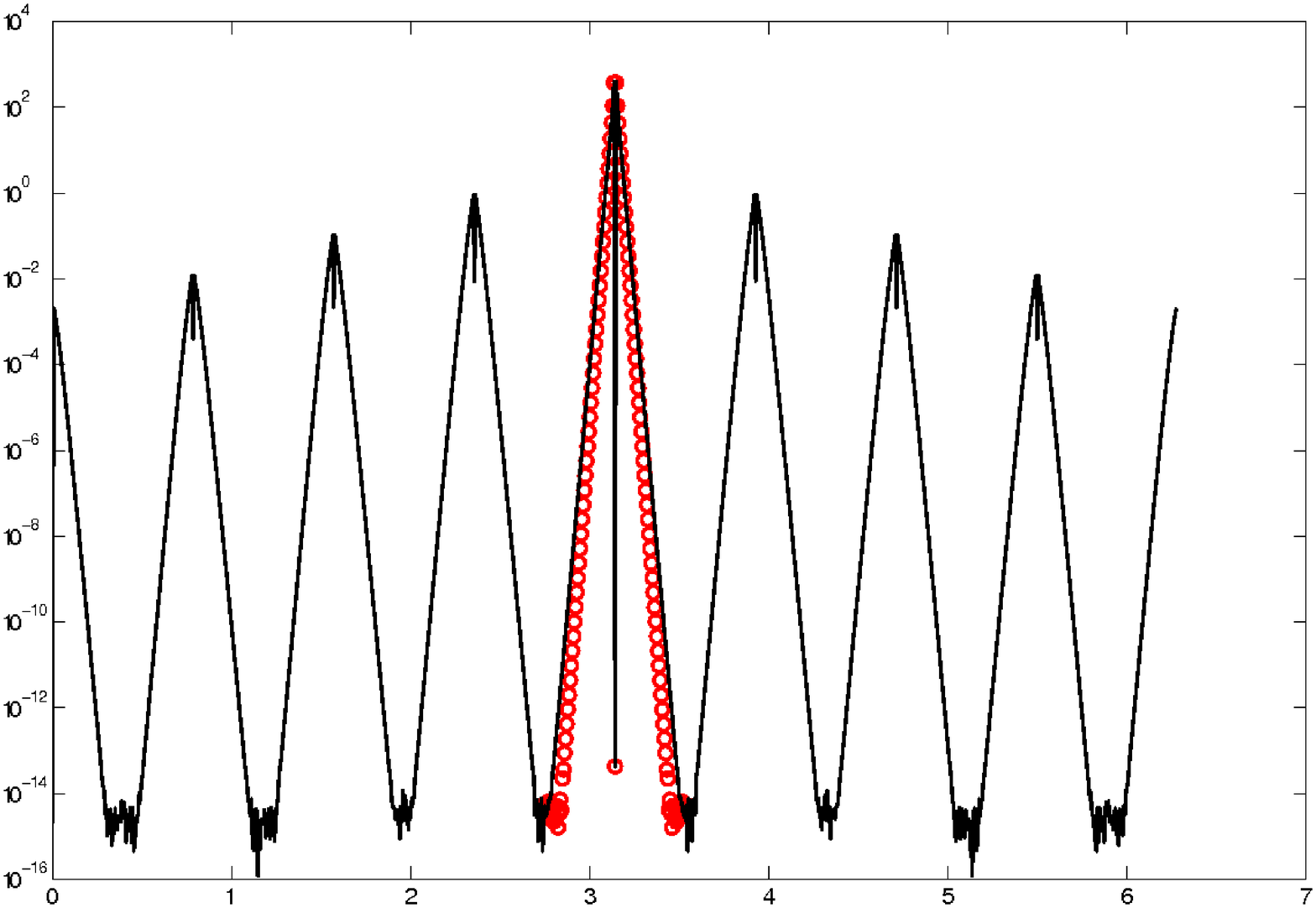}
}
\caption{Viscous Burgers equation with highly oscillatory diffusion. The solution is shown in $x$ space, zoomed in version in $x$ space, and in the $\widehat{u}$ space. The solutions are shown on a 1024 grid, with $dt=7.6e-6$, $dx=6.1e-3$, and $\lambda=6.3e-5$.}
\label{fig:Burgers}
  \end{center}
\end{figure*}

 For Figure \ref{fig:Burgers}, the diffusion coefficient is chosen as:
\begin{equation*}
a(x)=0.075 \exp\left({\frac{0.65+0.2 \cos(x)}{1+0.7 \sin(128x)}}\right)
\end{equation*}

\noindent The convolutions in the diffusion and non-linear terms are done in the spectral domain, rather than by other methods such as the psuedo-spectral method. The various dynamics can be seen in the spatial and in the spectral plots (see Figure \ref{fig:Burgers}). The true, sparse, and low frequency solutions are plotted in the space in Figure \ref{fig:Burgers}(a-b). The low frequency projection is done by thresholding any coefficients outside of a particular range. Specifically, the number of low wavenumbers retained is the same as the sparse solution, although their identities are dramatically different. The sparse solution captures the local and global behaviors of the solution more accurately than the low frequency projected solution.  In Figure \ref{fig:Burgers} (c-d), the spectrum of the true solution is compared to the sparse and low frequency spectra, respectively. The local peaks in the spectra are related to the wavenumbers in the diffusion coefficient $a(x)$ and the harmonics induced by the non-linear advection term. Notice that in this case, each of the distributions in the spectral domain do not decay as rapidly as in the parabolic case. The sparse solution contains 130 coefficients out of a total possible 1024, about $12.7\%$.

\subsection{Vorticity Equations}

The vorticity equation we consider is derived from the (2-d) incompressible Navier-Stokes equation \cite{majda2001vorticity}:

\begin{equation}
\partial_t u+ (\nabla^{\perp} \Delta^{-1} u) \ \cdot \nabla u=  \gamma \Delta u + f
\label{eq:Vort}
\end{equation}

\noindent where $u$ is the vorticity (not the velocity). Similarly to equation \ref{eq:Burgers}, equation \ref{eq:Vort} exhibits three separate phenomenon: 1. Smoothness from the diffusion term on the RHS, 2. Small scale oscillations induced from the high frequencies in the source term $f$, and 3. Non-linear interactions between frequencies from the advection term on the LHS. However, since the operator $\nabla^{\perp} \Delta^{-1}$ is smoothing (in some sense), the advection term can be viewed as less non-linear than the one found in equation \ref{eq:Burgers}. In terms of the numerical method, the operator $\nabla^{\perp} \Delta^{-1}$ dampens the coefficients by a factor which acts as $|k|^{-1}$.

For the numerical implementation, the diffusion term is discretized using Crank-Nicolson while the advection term is lagged. Since the operators are diagonalized in the coefficient basis, the steps can be invertible and lead to a simple updating scheme. 

\begin{eqnarray*}
  \widehat{v} ^{n+1}&=& \frac{2 dt}{2+\gamma dt |k|^2} \left( i k^{\perp} \left( |k|^{-2} u^{n} \right) \ast  i k \widehat{u}^{n} + \widehat{f} \right)\\
& & \ \ \ + \frac{2-\gamma dt |k|^2}{2+\gamma dt |k|^2} \  \widehat{u}^{n} 
\end{eqnarray*}

\begin{figure*}[tb]
\begin{center}
\subfigure[ True Solution in $x$ space]{
\includegraphics[width=2in]{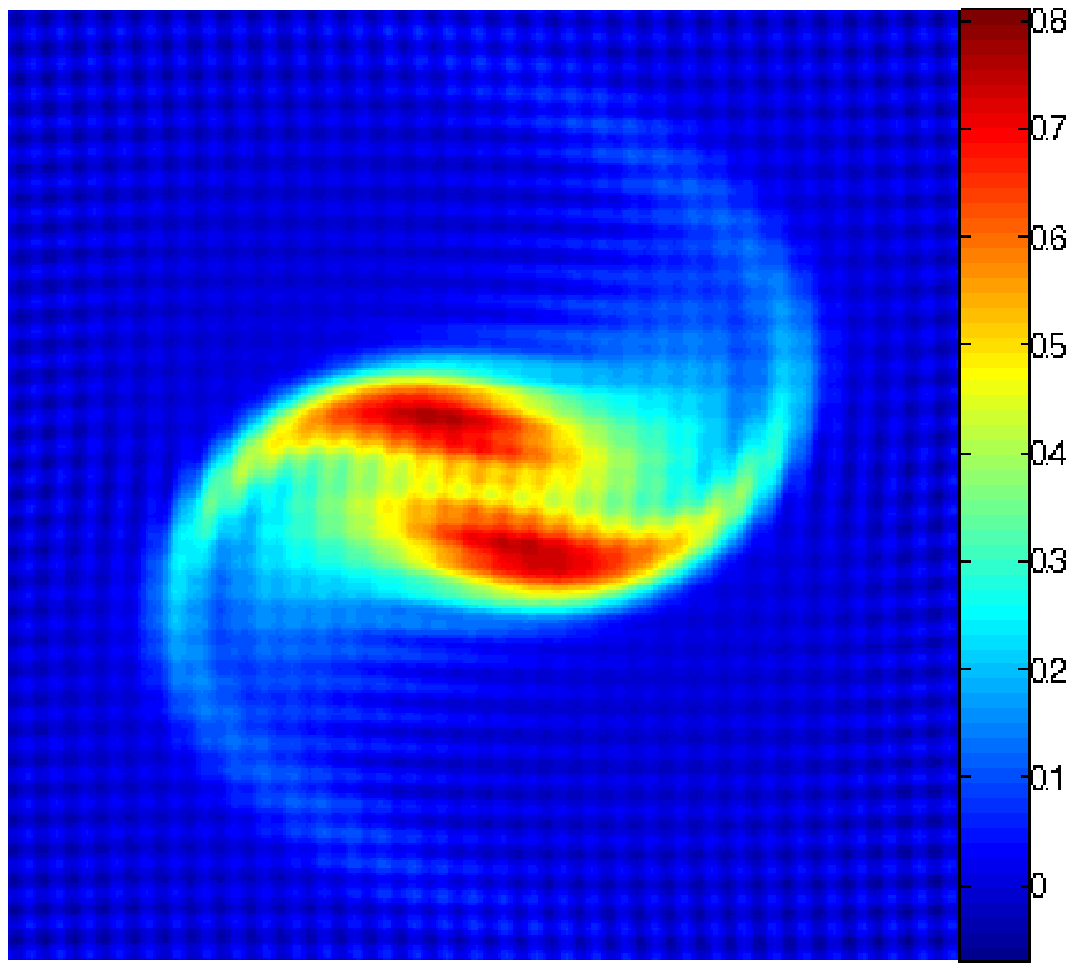}
} \hspace{0.2cm}
\subfigure[Sparse  Solution in $x$ space]{
\includegraphics[width=2in]{NS_sparse_256.eps}
}\hspace{0.2cm}
\subfigure[Low Frequency Solution in $x$ space]{
\includegraphics[width=2in]{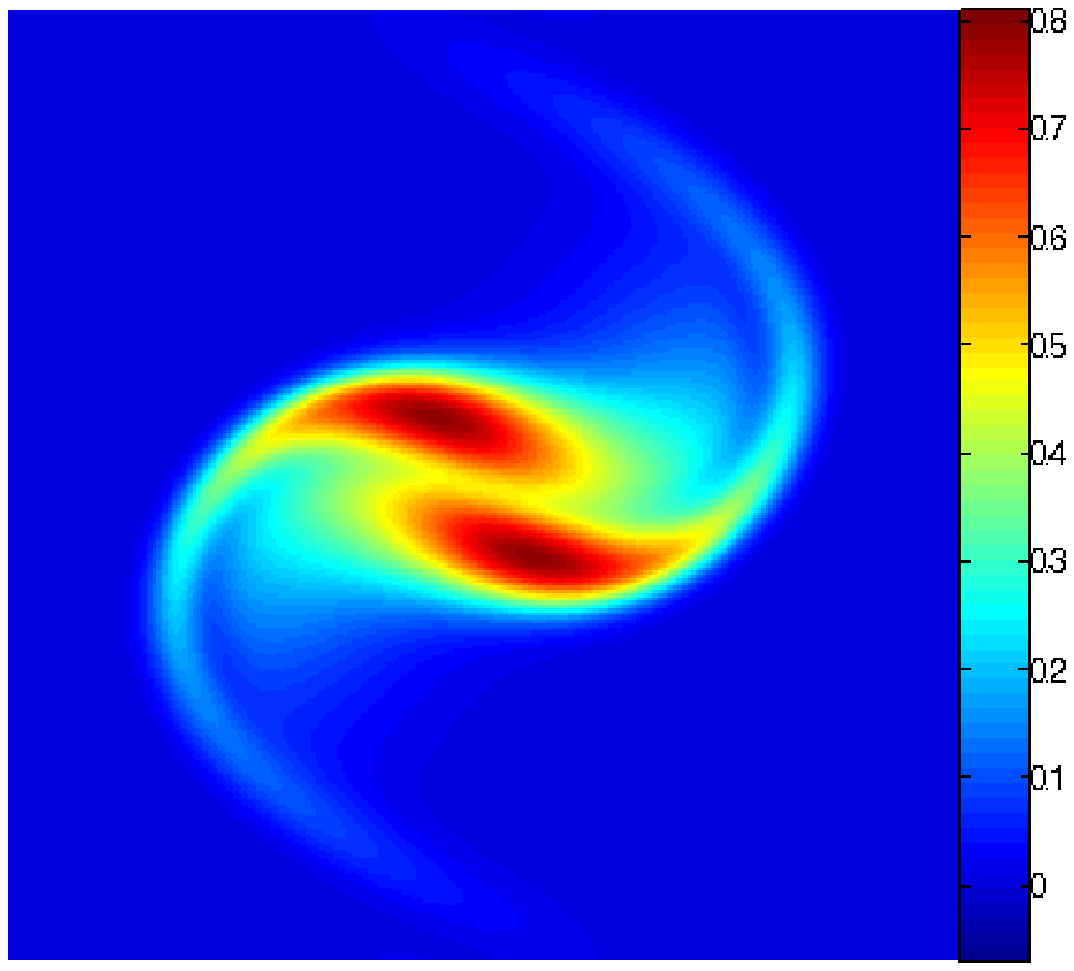}
}
\subfigure[True Solution in $\widehat{u}$ domain]{
\includegraphics[width=2in]{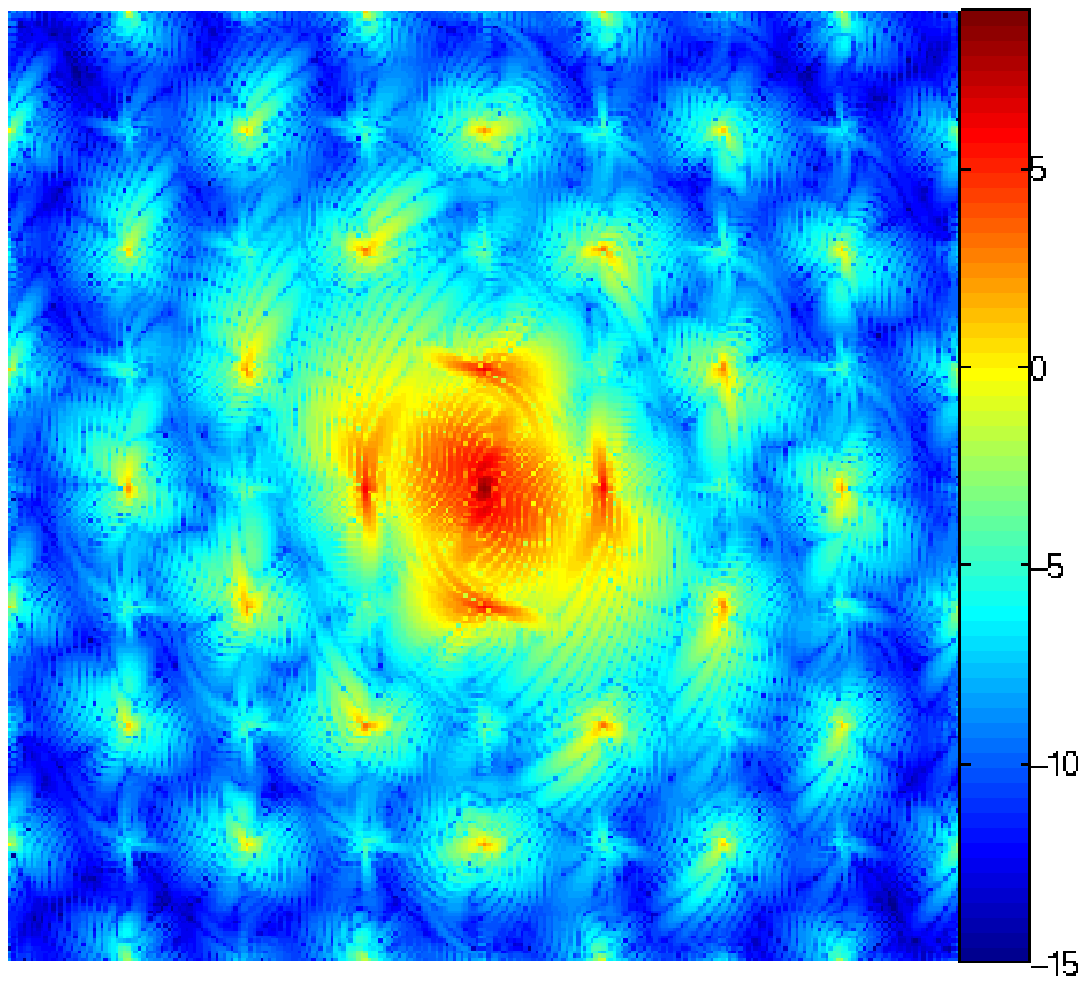}
}
\subfigure[Sparse Solution in $\widehat{u}$ space]{
\includegraphics[width=2in]{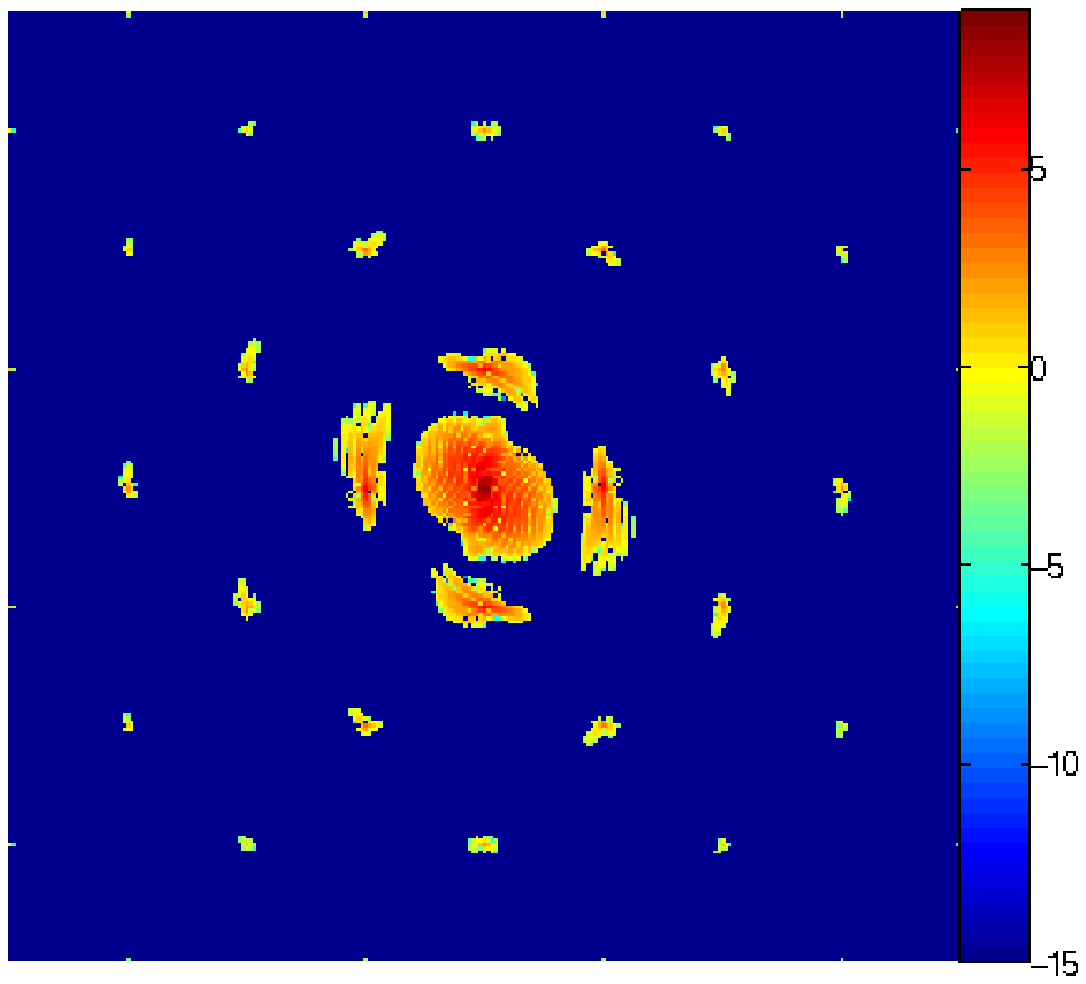}
}
\subfigure[$L^2$ (solid) and $L^\infty$ (dashed) Errors verse grid step size]{
\includegraphics[width=2.2in]{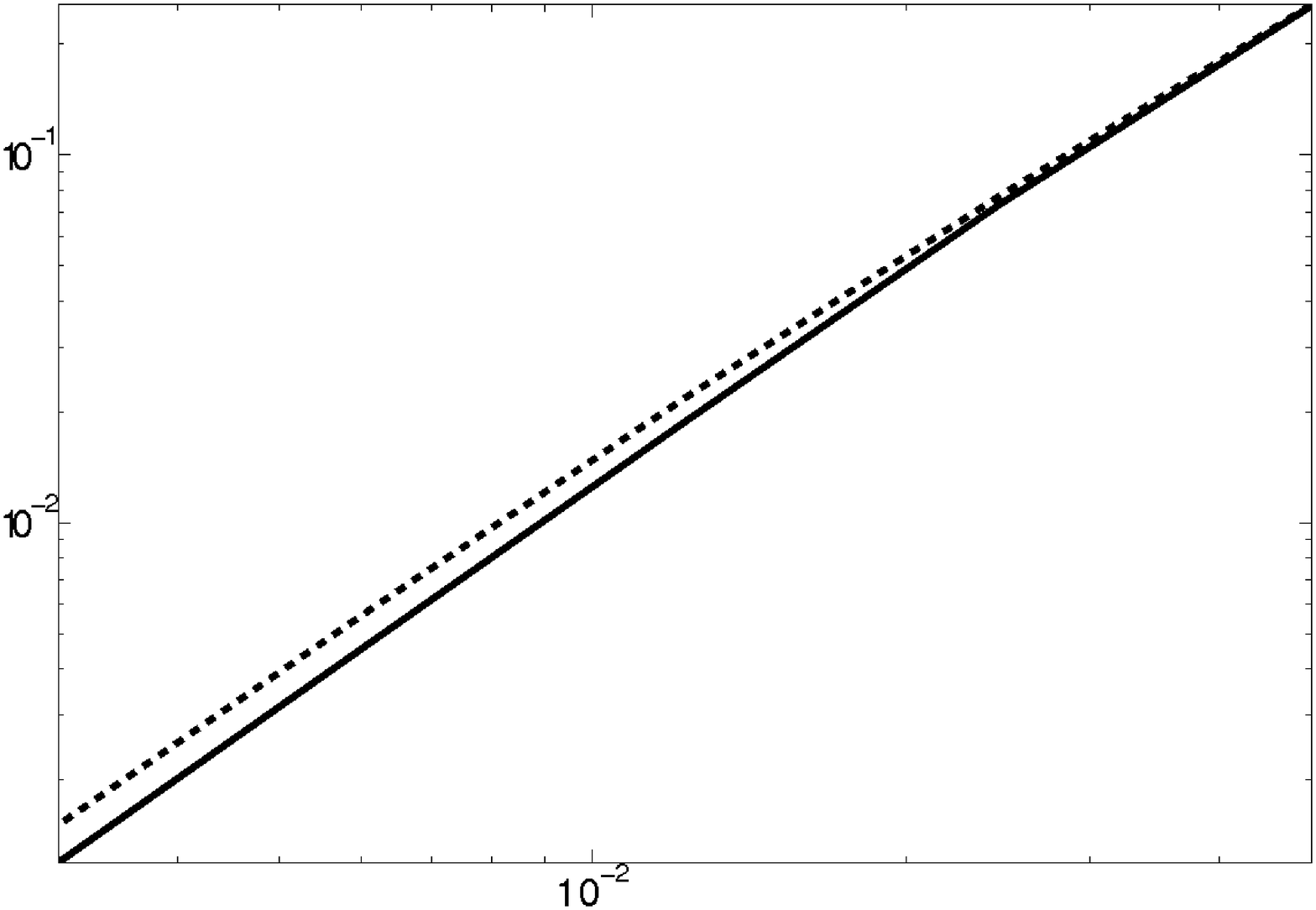}
}
\caption{Vorticity equations with high frequency source term. The solution is shown in $x$ space, zoomed in version in $x$ space, and in the $\widehat{u}$ space. The solutions are shown on a 256 by 256 grid, with $dt=0.025$, $dx=0.0245$, $\gamma=0.001$, and $\lambda=0.0497$.}
\label{fig:Vort}
  \end{center}
\end{figure*}

\noindent For Figure \ref{fig:Vort}, the forcing term is chosen to be:
\begin{equation*}
 f(x,y)=0.025 \frac{\sin(32x)+\sin(32y)}{1+0.25 \left( \cos(64x)+\cos(64y) \right)}
\end{equation*}

 \noindent The standard stability condition is used for choosing the time steps in order to insure capture of all small scale behaviors. In Figure \ref{fig:Vort} (a-b), the true and sparse solutions are plotted in the spatial domain. Notice the oscillations introduced by the source term interact with the two vortex patches and thus contribute to the global behavior of the solution. The spectra of the true and sparse solution is plotted in Figure \ref{fig:Vort} (c-d), where the low wavenumbers are located in the middle of the domain. The sparse solution retains about $3.95\%$ of the coefficients. In the sparse spectrum, the coefficients are located throughout the domain, including the highest frequency itself (seen on the boundary of the spectral domain). In Figure \ref{fig:Vort} (e), the $L^2$ and $L^\infty$ error are shown to decrease as the resolution increases. This sparse solutions, as well as the other examples presented here, converge as the spatial discretization goes to zero. 

It was observed that as the dimension increases, the sparsity of the solution also increases (since it is proportional to the product of the sparsities in each dimension). Thus the method scales well with dimension.

It is worth noting that in a related work, wavelet hard thresholding was used to separate coherent and incoherent structures in turbulent flows \cite{farge2001coherent}.

\setcounter{subfigure}{0}

\section{Discussion}
In the examples from the previous section, the PDEs contained a mixture of multiscale properties with a diffusion term. The combination of non-linear and oscillatory terms created a large range of wavenumbers in the solution, while the dynamics produced a range of amplitudes. This gives the necessary structure for sparsity with respect to the Fourier basis. In general, the highest order derivative will determine the appropriate basis in which the solutions could be sparse.

If the spectrum is more localized, \textit{i.e.} non-zero regions in the low frequency regime, then the proposed model can better condition the numerical method.  Empirically, the shrinkage operator acts as a non-linear filter on the coefficients.  It was observed that for a fixed $C$ and $p$, where $\lambda=C dt^p$, the numerical updates presented here with $dt$ larger than theoretically and numerically possible in the original scheme will converge. In the case of the vorticity example, $dt$ can be taken much larger when soft thresholding is applied than in the standard scheme. Also, the non-linear filter seems to reduce parasitic modes and spurious oscillations found in spectral approximations for linear  and for non-linear slightly viscous   hyperbolic equations (see Figures \ref{fig:Hyper} (b) and (c)) and  \ref{fig:Burgers} (g).

One key point  is that our method works by fully resolving the solution. Its efficiency is gained by omitting modes that are insignificant. This requires that $\lambda$ is small enough that the filter does not annihilate essential features.  For example, if the initial data is smaller than $\lambda$ for a particular unstable mode, then our approximate solution will not match the true dynamics. As the grid is refined, the mode will be captured (since $\lambda$ decreases as $\Delta x$ decreases).

In terms of complexity, each iteration is dominated by the convolution step. The convolution in the coefficient domain (spectral domain) can be done explicitly over the $n_s(t)$-sparse vectors rather than transforming back onto the spatial grid, which is $\mathcal{O} \left( n_s(t)^2 \right)$ at each step. When $n_s(t)^2 << N \mathrm{log} N$, convolving in the spectral domain rather than transforming back and forth between domains decreases the computational cost of the algorithm. Knowing \textit{a priori} the maximum sparsity, \textit{i.e.} $n_{s,max}=\max_t n_s(t)$, faster routines and transforms could be optimized for specific problems and applications.  For example, one can optimize the transform between the spatial and coefficient domains knowing the given sparsity at the current step and the non-trivial coefficients' identities. In the linear cases, as in equation \ref{eq:ParabolicUpdate}, the operation $\widehat{a} \ast$ can be stored as a large but sparse matrix, reducing the updates to a sparse matrix - sparse vector operation at every iteration. Our goal in this work is to formulate a PDE solver that promotes sparsity. In future work we will present a study of the computational complexity and speed.

When the dynamics are dominated by a linear term, for example high viscosity, the identities of the non-trivial coefficients settle over time. This was also observed in the non-linear cases, but over a longer time period. This would enable creation of a sparse basis for elliptic equations (e.g., with oscillatory coefficients).

In many of the cases here, it is possible to get hyper-sparse solutions (those with 1\% or fewer coefficients) at the cost of accuracy. 

\section{Conclusion}
In this work, we have proposed a method to fully resolve the solutions of multiscale PDEs while only retaining important modes. The reduced dynamics created by the sparse projection properly captures the true phenomena exhibited by the solution. The sparse projection amounts to a shrinkage of the coefficients of the updated solution at every time step. There exist many possibilities for using the sparsity to optimize individual algorithms and create faster, more efficient computational routines.








\section{Acknowledgments}
The research of H. Schaeffer was supported by the Department of Defense (DoD) through the National Defense Science and Engineering Graduate Fellowship (NDSEG) Program. The research of S. Osher was supported by ONR: N00014-11-1-719. The research of R. Caflisch was supported by DOE: DE-FG02-05ER25710. The research of C. Hauck is sponsored by the Office of Advanced
Scientific Computing Research; U.S. Department of Energy. The work was performed at the Oak Ridge National Laboratory, which is managed by
UT-Battelle, LLC under Contract No. De-AC05-00OR22725.  Accordingly, the
U.S. Government retains a non-exclusive, royalty-free license to publish or
reproduce the published form of this contribution, or allow others to do so,
for U.S. Government purposes.




\bibliographystyle{plain}
\bibliography{SparsePDE}










\end{document}